\documentclass{article}
\usepackage{amsfonts,amssymb,latexsym,amsmath}
\newcounter{num}[section]
\newcommand{\Num}{\refstepcounter{num}%
\textbf{\arabic{section}.\arabic{num}}}

\newcommand{\Theorem}{\textbf{Theorem~}}
\newcommand{\Proof}{\textbf{Proof}}
\newcommand{\Def}{\textbf{Definition~}}

\newcommand{\Lemma}{\textbf{Lemma~}}

\newcommand{\Prop}{\textbf{Proposition~}}
\newcommand{\Cor}{\textbf{Corollary~}}

\newcommand{\Ax}{{\mathfrak A}}
\newcommand{\Bx}{{\mathfrak B}}

\newcommand{\Kc}{{\cal K}}
\newcommand{\Oc}{{\cal O}}
\newcommand{\Xc}{{\cal X}}

\newcommand{\Ch}{{{\mathfrak C}{\mathfrak h}}}
\newcommand{\al}{{\alpha}}

\newcommand{\la}{{\lambda}}
\newcommand{\La}{{\Lambda}}
\newcommand{\tG}{{\widetilde{G}}}

\newcommand{\tGe}{{\widetilde{G_e}}}

\newcommand{\UTn}{{\mathrm{UT}}(n,\Fq)}
\newcommand{\Tn}{{\mathrm{T}}(n,\Fq)}
\newcommand{\Irr}{{\mathrm{Irr}}}

\newcommand{\Fq}{{\Bbb F}_q}
\newcommand{\Cb}{{\Bbb C}}

\newcommand{\ind}{{\mathrm{ind}}}

\newcommand{\xitl}{{\xi_{\theta, \la}}}
\newcommand{\xitlp}{{\xi_{\theta, \la'}}}
\newcommand{\xitla}{{\xi_{\theta, a\la}}}
\newcommand{\xitpl}{{\xi_{\theta', \la}}}

\newcommand{\chitl}{{\chi_{\theta, \la}}}
\newcommand{\chitlp}{{\chi_{\theta, \la'}}}
\newcommand{\chitla}{{\chi_{\theta, a\la}}}
\newcommand{\chitpl}{{\chi_{\theta', \la}}}
\newcommand{\chitplp}{{\chi_{\theta', \la'}}}

\newcommand{\Jlr}{{J_{\la,\mathrm{right}}}}
\newcommand{\Jlpr}{{J_{\la',\mathrm{right}}}}

\renewcommand{\leq}{\leqslant}

\newcommand{\rt}{{\mathrm{right}}}
\newcommand{\lt}{{\mathrm{left}}}

\begin{document}
\Large

\title{Supercharacters for the finite groups of  triangular type}
\author{A.N.Panov
\thanks{The paper is supported  by RFBR grant  14-01-97017-Volga Region-a}\\
Samara State University\\ apanov@list.ru}
\date{}
 \maketitle

{\small {\bf \hspace{0.3cm} MSC:} 20C05, 22E27}

{\small {\bf \hspace{0.3cm} Keywords:} group representations, supercharacter theory, triangular group, orbit method

\abstract{
 We construct the supercharacter theory for the finite groups of triangular type. Its  special case is the supercharacter theory for  algebra groups of P.Diaconis and I.M.Isaacs.  The supercharacter analog of the A.A. Kirillov formula for irreducible characters is proved.
 }

\section{Introduction}

 It is well known  ~\cite{Dr} that  for some unipotent groups, like  the unitriangular group $\UTn$, the problem of classification of irreducible representations is an extremely difficult  "wild"\, problem.
 In the series of  papers ~\cite{A1,A2,A3,A4}, C.A.M.Andre  constructed the theory of so called basic characters of the group  $\UTn$ (see also \cite{Yan, P1}). These characters are not irreducible in general, nevertheless they  have many common features with the theory of irreducible characters.  Later on this approach was developed by  P.Diaconis and  I.M.Isaacs in the paper   ~\cite{DI}. They suggested to change the problem of classification of irreducible representations by the problem of construction so called theory of supercharacters.  Each group admits  different superchracter theories including the trivial ones.   More precisely, the problem is to construct  a supercharacter theory which as much as possible close to  the theory of irreducible  representations.  The theory of basic characters  is an example of a supercharacter theory; up today nobody proposed the better one for   $\UTn$.
 In mentioned above paper  ~\cite{DI}, P.Diaconis and  I.M.Isaacs constructed  the supercharacter theory for the algebra groups  (by definition, an algebra group is a group of the form  $1+J$, where  $J$ is an associative finite  dimensional nilpotent algebra over the finite field   $\Fq$), which in the special case $\UTn$ coincides with  the theory of basic characters.

 In the present paper, we consider  the class of finite groups of triangular type. Algebra groups are the special case of these groups. In the sections 1-5 of the present paper, we construct the supercharacter theory for the groups of this type  (see Theorem \ref{main}), which in the case of algebra groups coincides  the theory of P.Diaconis and  I.M.Isaacs.
 An other special case is the supercharacter theory for the groups of invertible elements of the reduced algebras over a  finite field  (see  Example 3) that was constructed in the author's paper  ~\cite{P1}. Notice that the present supercharacter theory affords much better approximation of the theory of irreducible representations than the general  construction of the paper  \cite{H}.
 In the final  section  6, we construct the supercharacter analog of the A.A.Kirillov formula for irreducible characters of unipotent groups (see Theorem \ref{kirillov}).

Following the paper \cite{DI}, we recall the main definitions and statements of the supercharacter theory. Let $G$  be a finite group. Notice that we  denote the unit element  by $1$. Consider the system of characters (representations)  $$\Ch = \{\chi_\al|~ \al\in \Ax\} $$  of the group  $G$. \\
 \Def\Num. The system of characters  $\Ch$  defines a supercharacter theory of the group  $G$
 if the characters of  $\Ch$ are pairwise disjoint and there exists a partition  $$ \Kc = \{K_\beta|~ \beta\in\Bx\}$$ of the group $G$ such that :\\
S1) $|\Ax|= |\Bx|$;\\
S2) ~ each character  $\chi_\al$ is constant on each $K_\beta$;\\
S3) ~ $\{1\} \in \Kc$.
\\
We refer to each character  of  $\Ch$ as  a {\it supercharacter} and to subset of  $\Kc$ as  a {\it superclass}.

 Denote by  $\Xc$ the system of subsets  $\{X_\al\}$ of  $\Irr(G)$, where each $X_\al$  is  the subset of all irreducible components of  $\chi_\al$. One can understand the importance of the condition  S3) by  the following statement. \\
\Prop\Num\label{lll} \cite[Лемма 2.1]{DI}. Suppose that  the system of disjoint characters  $\Ch$ and the partition  $\Kc$ subject conditions  S1) and  S2). Then
the condition  S3) is equivalent to the following condition \\
S4)~ the system of  subsets $\Xc$ form a partition of $\Irr(G)$ and each character  $\chi_\al$ is up to a constant multiple equals to the character  $$\sigma_\al = \sum_{\psi\in X_\al} \psi(1)\psi. $$
\Cor\Num. The system of supercharacters is uniquely determined up to constant factors by the partition  $\Irr(G) = \bigcup X_\al$.

The lemma  \ref{lll} afforded the new definition:  a supercharacter theory of the  group $G$ is defined by two partitions  $\Irr(G) = \bigcup X_\al$ and  $G=\bigcup K_\beta$ such that each  $\sigma_\al$  is constant  on each  $K_\beta$. This definition is usually used in the general theory. If one desired to construct a special supercharacter theory it is better to use the first definition since,  as we say above, the problem of classification of all irreducible characters may happens to be a "wild" problem.

Notice also some important properties of the systems  $\Ch$ and $\Kc$.\\
\Prop\Num ~\cite[Теорема  2.2.]{DI}. \\
1) Each superclass is a union of the classes of conjugate elements.\\
2) The partition  $\Kc$ is uniquely defined by the partition $\Xc$ and vice versa.\\
3) The principal character is a supercharacter (up to a constant factor).

\section{Finite groups of triangular type}

 Let $H$ be a group and $J$ be an associative algebra over a field  $k$. There defined the left  $h,x\to hx $ and  right $h,x\to xh$ actions of  $H$  on $J$, and these actions commute. Suppose that the following conditions hold:
\begin{enumerate}
\item $h(xy)=(hx)y$ and $(xy)h=x(yh)$,\\
\item $x(hy)=(xh)y$.
\end{enumerate}
The set  $$G=H+J=\{h+x:~ h\in H, ~ x\in J\}$$
is equipped  by the associative operation of multiplication:
\begin{equation}\label{mult}
g_1g_2=(h_1+x_1)(h_2+x_2)=h_1h_2+ h_1x_2+x_1h_2+x_1x_2.
\end{equation}
If  $J$ is a nilpotent algebra over the field  $k$, then  $G$ is a group with respect to the operation (\ref{mult}). If the group  $H$ is finite, and $k$ is a finite field, and  $J$ is a finite dimensional  nilpotent algebra over the field $k$, then   the group $G$ is  finite. \\
\Def\Num. We refer to the constructed group   $G$ as a \emph{finite group of triangular type} if    $H$ is  abelian   and  $\mathrm{char}\, k$ does not divide  $|H|$.

Let  $G=H+J$ be a  finite group of triangular type. According to Maschke's theorem, the group algebra  $kH$  is commutative and semisimple,  and  therefore is a sum of fields.
There exists a system of primitive idempotents  $\{e_1,\ldots, e_n\}$ such that   \begin{equation}\label{gralg} kH=k_1e_1\oplus\ldots \oplus k_ne_n, \end{equation} where $k_1,\ldots,k_n$ are  field extensions of  $k$.

The direct sum  $A=kH\oplus J$ has an algebra structure  with respect to multiplication  (\ref{mult}). The group  $G$ is a subgroup of the group  $A^*$ of invertible elements of $A$, which is considered in the example 3.
Notice that the group  $G$  decomposes into a product $G=HN$ of the subgroup  $H$ and the normal subgroup   $N=1+J$ that is an algebra group.

Consider the group  $\tG$ of the triples  $\tau=(t,a,b)$, where  $t\in H$,~ $a,b\in N$,  with operation  $$(t_1,a_1,b_1)\cdot (t_2,a_2,b_2) = (t_1t_2,~ t_2^{-1}a_1t_2a_2, ~ t_2^{-1}b_1t_2b_2).$$
Define a representation of the group $\tG$ in  $J$  by the formula $$ \rho(\tau)(x) = taxb^{-1}t^{-1}.$$
The representation of $\tG$ in the dual space  $J^*$ is defined as usual  $$\rho^*(\tau)\la(x) = \la(\rho(\tau^{-1})(x)).$$
 There are also the  left $b\la(x)= \la(xb)$ and  right $\la a(x) = \la(ax)$ linear actions  of $G$ on   $J^*$. Then  $\rho(\tau)(\la) =
tb\la  a^{-1}t^{-1}$.

  For every idempotent  $e\in A$, denote by  $A_e$ the subalgebra $eAe$. The subalgebra  $J_e=eJe$ is the radical in  $A_e$. Denote $e'=1-e$. The Pierce decomposition takes place
 $$J= eJe\oplus eJe'\oplus e'Je\oplus e'Je'.$$  We identify  the dual space $J_e^*$ with the subspace in $J^*$ of all linear forms equal to zero on all components of the Pierce decomposition apart of the first one. \\
\Def\Num. \\
1) We say that  an element    $x\in J$
  is singular if $x\in J_e$ for some idempotent  $e\in A$, ~ $e\ne 1$.
   Otherwise  $x$ is   regular.  \\
2) We say that  an element  $\la\in J^*$
is singular if $\la\in J_e^*$ for some idempotent  $e\in A$, ~ $e\ne 1$. Otherwise $\la$ is  regular.\\
\Lemma\Num\label{ccc} \cite[Lemma 2.2]{P1}. The following conditions are equivalent:
1) the element  $x\in J$ (respectively,  $\la\in J^*)$ is  singular;
2) there exists $c\in A\setminus J$ such that  $cx = xc = 0$ (respectively, $c\la = \la c = 0$).
\\
\Prop\Num. If an element $x\in J$ (respectively, $\la\in J^*$) is  singular,  then
each element of its   $\rho(\tG)$-orbit (respectively, $\rho^*(\tG)$-orbit) is  singular. Similarly for  regular elements.\\
\Proof. We prove for  $x\in J$. The case  $\la\in J^*$ is similar.

Let $x$ be a  singular element. There exists an idempotent  $f\in A$, ~ $f\ne 0$ such that  $fx=xf=0$. It suffices to prove that  the element  $ax$ (respectively, $xa$) is  singular for every  $a\in N$. For  $c=fa^{-1}$ we have $c\notin J$ and $c(ax)=(ax)c=0$. The lemma  \ref{ccc}
implies that $ax$ is  singular. $\Box$\\
\Cor\Num\label{einS}. For any  singular $\tG$-orbit $\Oc\subset J$ (respectively, $\Oc^*\subset J^*$) there exists an idempotent  $e\in kH$ such that  $\Oc\cap J_e\ne \varnothing$ (respectively, $\Oc^*\cap J_e^*\ne \varnothing$).\\
\Proof. For any idempotent  $e\in A$ there exists  $g\in A^*$ such that $geg^{-1}\in kH$ \cite[теорема 4.1]{DK}. Since $g=sa$, where $s\in (kH)^*$ и $a\in N$, we have $aea^{-1}\in kH$.  $\Box$

Recall that the subgroup   $H$ is abelian; then  $he=eh=ehe$ for any $h\in H$. The subset  $H_e=eHe$ is a subgroup in the group of invertible elements of the algebra  $A_e$. The group $G_e=eGe=H_e+J_e$ is a group of triangular type, and it is associated with the algebra  $A_e$ in the same way as    $G$ is associated with  $A$.
Similarly to the group $\tG$, we define the group  $\tGe$. Notice that the map  $h\to he$ is a homomorphism of the group  $H$ onto  $H_e$;  its kernel is the subgroup
\begin{equation}\label{he} H(e) =\{h\in H: ~ he=e\}\end{equation}
 Let  $H_{y,\lt}$, ~  $H_{y,\rt}$,~ $H_{\la,\lt}$, ~  $H_{\la,\rt}$  be stabilizers of the elements $y\in J$ and $\la\in J^*$ with respect to the left and right actions of  $H$.\\
\Lemma\Num\label{leftright}.\\
1)~  $H(e)   = H_{y,\mathrm{right}}\cap H_{y,\mathrm{left}}$ for any regular  $y\in J_e$.\\
2)~ $H(e) = H_{\la,\mathrm{right}}\cap H_{\la,\mathrm{left}}$ for any regular  $\la\in J_e^*$.\\
\Proof. Let us prove  2), the statement 1) is similar.  If  $h\in H(e)$, then $he=eh=e$. For any element  $\la\in J_e^*$, we obtain
$h\la=he\la=e\la=\la$ and $\la h=\la eh=\la e=\la$. Then   $H(e)$ is contained in the intersection of the left and right stabilizers.

Let us prove the contrary.  Let $h\la = \la h=\la$.  As $\la\in J_e^*$, we have
$e\la=\la e=\la$. Then for  $c=he-e$ from  $A_e$, we obtain  $c\la = \la c =0$. If $c\ne 0$, then $c\in A_e\setminus J_e$;  the lemma  \ref{ccc} implies that  $\la$ is singular in $J_e^*$. This contradicts to choice of  $\la$.
Therefore, $c=0$ and $he=e$. That is  $h\in H(e)$. $\Box$
\\
\Lemma\Num\label{Oc}. Let  $e,f$ be two idempotents in  $kH$.  Then\\
1)~ if the  $\tG$-orbit $\Oc$ has nonempty intersection with both   $J_e$  and $J_f$, then the orbit $\Oc$
has nonempty intersection with $J_{ef}$; \\
2) two elements   $x $ and  $y$ in  $J_e$  belong to a common  $\tG$-orbit $\Oc$ if and only of  $x$ and  $y$ belong to a common $\tGe$-орбите. \\
\Proof.  \\
1) Let $x\in \Oc\cap J_e$, ~$y\in \Oc\cap J_f $ (i.e. $ex=xe=x$ and $fy=yf=y$, and $y=haxbh^{-1}$, where $h\in H$ and $a,b\in N$). Consider the element   $z=eye$. Since $ef=fe$, we have $(ef)z=z(ef)=z$ and  $z=eye= ehaxbh^{-1}e = heaxbeh^{-1}$.
  The element $a_1= ea+e'$, where $e'=1-e$, obeys   $a_1x=eax+(1-e)x=eax$. As  $a_1 = 1\bmod J$, we have $a_1\in N$.
Similarly,   $b_1=be+ e'$ belongs to  $N$ and $xb_1=xbe+x(1-e) = xbe$. Then $z=ha_1xb_1h^{-1}$, where  $a_1,~b_1\in N$; that is
$z$ belongs to the  orbit   $\Oc$.\\
2a) Let  $x,y\in J_e$ ( i.e.  $ex=xe=x$ and $ey=ye=y$).
 Suppose that  $x,y$ belong to a common  $\tGe$-orbit. Then
 $y=(ehe)(eae)x(ebe)(eh^{-1}e)=h(eae+e')(exe)(ebe+b')h$. As we see above, the elements   $a_1=eae+e'$ and  $b_1=ebe+e'$ belong to $N$. We obtain $y=ha_1(exe)b_1h^{-1}=ha_1xb_1h^{-1}$. The elements $x,y$ belong to a common  $\tG$-orbit.\\
 2b) Let   $x,y$ of  $J_e$  belong to a common $\tG$-orbit,  i.e.   $y=haxbh^{-1}$ for some $a,b\in N$ and  $h \in H$. Then
$$y=eye=ehaxbh^{-1}e = (he)a(exe)beh^{-1} = (he)(eae)x(ebe)(h^{-1}e).$$ The element  $he$ belongs to  $H_e$ and its inverse element is  $h^{-1}e$. The element  $eae$ and $ebe$ are lying in  $N_e$. This proves that  $x$ and  $y$ belong to a common  $\tGe$-orbit. $\Box$\\
\Cor\Num\label{eee}. For any   $x\in J$ there exists a unique idempotent  $e\in kH$
such that  $\Oc(x)\cap J_e\ne \varnothing$ and    $\Oc(x)\cap J_e$ is a regular  $\tGe$-orbit in $J_e$. \\
\Proof. The statement follows from  1) and 2) of the previous lemma.
$\Box$

For following properties of orbits in  $J^*$ can be proved similarly.\\
\Lemma\Num\label{Oc*}. Let  $e,f$ be two idempotents  in  $kH$. Then \\
1)~ if the   $\tG$-orbit in  $\Oc^*$ has nonempty intersection with both $J_e^*$ and $J_f^*$, then $\Oc^*$
has nonempty intersection with $J_{ef}^*$; \\
2) two elements  $\la $ and  $\mu$ in $J_e^*$ belong to a common $\tG$-orbit $\Oc^*$ if and only if   $\la$ and $\mu$ belong to a common $\tGe$-orbit.\\
\Cor\Num. For any  singular element  $\la\in J^*$ there exists a unique edempotent  $e\in kH$ such that  $\Oc(\la)\cap J^*_e\ne \varnothing$, and   $\Oc(\la)\cap J^*_e$ is a  regular  $\tGe$-orbit in $J^*_e$.\\
\Cor\Num\label{cap}.~ 1) $\Oc(J_e)\cap \Oc(J_f) = \Oc(J_{ef})$,~~ 2) ~$\Oc(J_e^*)\cap \Oc(J_f^*) = \Oc(J_{ef}^*)$.
\Prop\Num\label{ereg}. The number of all regular  $\tG$-orbits  in  $J$ coincides with the number  of all  regular  $\tG$-orbits  in  $J^*$ (i.e.  $n_E(J) = n_E(J^*)$).\\
\Proof.
  Let $\{e_1,\ldots, e_n\}$ be a system of all primitive idempotents of the algebra  $S$. Denote $e'_i = 1-e_i$.
For each idempotent  $f\in S$, we denote by $l(f)$ the number of factors in the decomposition  $f=e'_{i_1}\cdots e'_{i_l}$. For  $f=1$, we take  $l(f)=0$.
The Corollary   \ref{cap} implies that
$\Oc(J_{f})$  coincides with the intersection  $\bigcap \Oc(J_{\phi})$, where $\phi $ runs through the set of idempotents   $\{e'_{i_1},\ldots, e'_{i_l}\}$.

The set of all singular orbits  in   $J$ is a union of all $\Oc(J_{e'_i})$, where $1\leq i\leq n$. Hence
$$ n_E(J) = \bigoplus (-1)^{l(f)} n(J_f),$$
where  $f$ runs through the set of all idempotents of  $S$. The similar formula is true for the orbits in  $J^*$:
$$ n_E(J^*) = \bigoplus (-1)^{l(f)} n(J^*_f), $$
where $f$ is also runs through the set of all idempotents of $S$.
For any finite subgroup of linear operators in a linear space $V$ over a finite field, the number of orbits in  $V$ coincides with the number of orbits in the dual space  $V^*$ (см.\cite[Lemma 4.1]{DI}). For the orbits of the group  $\tilde{G_f}$ in  $J_f$  and $J^*_f$, we conclude that $n(J_f) = n(J^*_f)$. Therefore  $n_E(J)=n_E(J^*)$. $\Box$

\section{Superclasses}

For every   $g\in G$ and $\tau = (t, a, b)\in \tG$, where  $ t\in H$ and $a,b\in N$, we define an element  $R_\tau(g)\in A$ by the formula
\begin{equation}\label{RRR}
    R_\tau(g) = 1+ ta(g-1)b^{-1}t^{-1}
\end{equation}
 If $g=h+x$, ~ $h\in H$,~ $x\in J$, then $R_\tau(g)= h\bmod J$. Hence  $R_\tau(g)\in G$. The formula  (\ref{RRR}) defines an action of the group  $\tG$ in $G$.
 {\it A superclass }  $K(g)$ is a $\tG$-orbit of the element  $g\in G$.

Two elements of $kH$ are   {\it associated } if they differs by a multiple from $(kH)^*$. By  (\ref{gralg}), we see that each element  of  $kH$  is associated with some idempotent of  $kH$.
\\
\Theorem\Num\label{superclass}. Let  $g=h+x$, where $h\in H$ and $x\in J$. Let  $f'$ be the idempotent in $kH$ associated with  $s=h-1$, and  $f=1-f'$.  Then\\
1) there exists the element $h+y$ in $K(g)$ such that  $y\in J_f$ (equivalent to $hy=yh=y$);\\
2) the elements  $h+y$ and $h+y'$, where $y, y'\in J_{f}$, belong to a common superclass if and only if  $y$ and $y'$ belong to a common $\tG_{f}$-orbit in  $J_{f}$.\\
\Proof.\\
{\bf Item 1}.  Let $g=h+x$ as above. To prove 1) it suffices to show that there exist  elements  $a,b\in N$ and $y\in J$ such that
$s+y=a(s+x)b$ and  $yf' = f'y = 0$. This reduces to the proof of the proposition that  for any  $x\in J$ there exist
$a,b\in N$  and  $y\in J$ such that  $f'+y=a(f'+x)b$ and $yf'=f'y=0$.\\
i) Let us prove that there exist  $u\in J$ and $a\in N$ such that  $f'+u=a(f'+x)$ and $uf'=0$.
 Take  $a=(1+x)^{-1}\in N$. Then
 $$
 uf'=\left((1+x)^{-1}(f'+x)-f'\right)f'= (1+x)^{-1}(f'+xf')-f' = 0.$$
 ii) Let  $u$ be as in i). Let us show that there exist $y\in J$ and $b\in N$ such that $f'+y = (f'+u)b$,~ $ b\in N$, and  $yf'=f'y=0$.  Take    $b=(1-f'u)\in N$. Then          $$y=(f'+u)(1-f'u)-f'=(1-f')u.$$ Hence $f'y=yf'=0$.\\
 {\bf Item 2}. We shall prove 2).
Easy to see that if $y$ and  $y'$ belong to a common  $\tG_{f}$-orbit, then  $h+y$ and  $h+y'$ belong to a common superclass. Let us prove the contrary.
Let  $h+y$ and  $h+y'$ belong to a common superclass.  Then there exist  $t\in H$, ~ $u,v\in J$ such that
\begin{equation}\label{hhh}
 h-1 + y' = t(1+u)(h-1 + y)(1+v)^{-1}t^{-1}.
\end{equation}
This implies
\begin{equation}\label{huv}
 (h-1 + y')t(1+v) = t(1+u)(h-1 + y)
\end{equation}
The element defined by the equality  (\ref{huv}) belongs to
$(h-1)t+J$. Subtract the element  $(h-1)t$ from the left anf right sides of (\ref{huv})  and multiply both sides by  $f$. Taking into account $f(h-1)=(h-1)f=0$, we obtain
$$fy't(1+v)f = ft(1+u)yf.$$
Since $fy=yf=y$ and $fy'=y'f=y'$,  we finally have
$$ y'ft(1+fvf) =  ft(1+fuf)y.$$
That is  $y$ and  $y'$ belong to a common $\tG_{f}$-orbit in  $J_{f}$. $\Box$\\
\Cor\Num\label{eeque}. Let $g=h+x$, ~$f$ and $f'$  be as in the theorem. Then \\
1)  there exists an idempotent  $e<f$ and regular element $y\in J_e$ (as an element in  $A_e$), such that  $h+y$ belongs to the same superclass  $K(g)$; denote by $\omega$ the orbit  of   $y$ with respect to  $\tG_e$;\\
2)~ the triple  $\beta = (e,h, \omega)$ is uniquely determines by the superclass  $K(g)$.\\
 \Proof.  For any  $\tG_{f}$-orbit  $\Omega$ in  $J_{f}$ there exists a unique idempotent  $e<f$ such that  $\Omega\cap J_e$ is a  regular  $\tG_e$-orbit in  $J_e$ (see corollary  \ref{eee}).
We can consider that the element  $y$ from the  Theorem \ref{superclass} to be a regular element  of $J_e$. This proves the statement   1).
The Theorem  \ref{superclass} and Lemma  \ref{Oc} imply 2). $\Box$

Denote by  $\Bx$ the set of triples  $\beta = (e,h, \omega)$, where
 $e$  is an idempotent in   $kH$, ~ $h\in H(e)$ (i.e. $he=e$),  and $\omega$ is a regular $\tG_e$-orbit in $J_e$. Choose $y\in\omega$. By Lemma  \ref{leftright}, the condition $h\in H(e)$ is equivalent to  $hy=yh=y$.  Denote by $K_\beta$ the superclass of  $h+y$. By Corollary \ref{eeque},  $K_\beta$  does not depend on the choice of $y\in \omega$.    \\
\Theorem\Num\label{supcl}. The  correspondence  $\beta\to K_\beta$ is a bijection of the set  of triples  $\Bx$ onto the set of all superclasses in $G$. The number of superclasses equals to
\begin{equation}
|\Bx| = \sum n_E(J_e)+|H(e)|,
\end{equation}
where  $e$ runs through the set of all idempotents of $kH$  for whom   the set  of regular orbits in $J_e$ is nonempty (i.e. $ n_E(J_e)\ne 0$).

\section{Supercharacters}

Consider the set   $\Ax$ of all triples $\al = (e,\theta, \omega^*)$, where
 $e$ is an idempotent in   $kH$, ~ $\theta $ is a linear character (one dimensional representation) of the subgroup  $ H(e)$,  and  $\omega^*$ is a  regular $\tG_e$-orbit in  $J^*_e$.
 Since the subgroup  $H(e)$ is abelian, the number of its linear characters is equal to the number of its elements. By proposition  \ref{ereg}, we have $|\Ax|=|\Bx|$.

 Let $\al = (e,\theta, \omega^*)\in \Ax$ and $\la\in \omega^*$.
 Denote
$$J_{\la,\mathrm{right}} = \{x\in J|~ \la x =0\},\quad\quad  N_{\la,\mathrm{right}} = \{ a\in G|~ \la a =\la\}.$$
Notice that  $N_{\la,\mathrm{right}} = 1+ J_{\la,\mathrm{right}}$.
Consider the subgroup  $G_\la = H(e)\cdot N_{\la,\mathrm{right}}$. The subgroup  $G_\la$ is a semidirect product of  $H(e)$ and the normal subgroup  $N_{\la,\mathrm{right}}$.  Any element  $g\in G_\la$ is uniquely presented in the form  $g=h+x$, where  $h\in H(e)$ and $x\in J_{\la,\mathrm{right}}$.

Let us fix  a nontrivial character  $t\to  \varepsilon^{t}$  of the additive group of the field $\Fq$ with the values in the multiplicative group  $\Cb^*$. By the triple  $\al = (e,\theta, \omega^*)$ and $\la \in\omega^*$, we define the linear character  of the subgroup  $G_\la$ as follows
\begin{equation}\label{}
    \xi_{\theta,\la}(g) = \theta(h)\varepsilon^{\la(x)},
\end{equation}
where $g=h+x$,~ $h\in H(e)$ and $x\in J_{\la,\mathrm{right}}$.
Let us show that $\xi=\xi_{\theta,\la}$ is really a linear character:
$$\xi(gg') = \xi((h+x)(h'+x')) = \xi(hh' + h'x + x'h  + xx') =$$
$$\theta(hh')\varepsilon^{\la(h'x)} \varepsilon^{\la(x'h)}\varepsilon^{\la(xx')} = \theta(h)\theta(h')\varepsilon^{\la(x)} \varepsilon^{\la(x')} = \xi(g)\xi(g').$$
We shall refer to the following induced character as a  {\it supercharacter}
\begin{equation}\label{indchi}
\chi_{\theta,\la} = \ind(\xi_{\theta,\la}, G_\la, G).
\end{equation}
\Prop\Num\label{supersuper}. The supercharacter   $\chi_{\theta,\la}$ is constant on superclasses.\\
\Proof. \\
{\bf Item 1}. Let $g\in G_\la$, ~$a\in N$. Let us show that   $g'=1+(g-1)a\in G_\la$ and $\chitl(g')=\chitl(g)$.

If $g\in G_\la$, then $g=h+x$, where  $h\in H(e)$ and  $x\in \Jlr$. Let  $a=1+u$ where  $u\in J$.
 Then $g'= h+y$, where  $h\in H(e)$ and $y=(h-1)u + x+ xu\in J$.  Since  $$\la h =\la,\quad \la(h-1)u=0, \quad \la x=0~~\mbox{and}~~ \la xu=0,$$  we obtain
$\la y = \la (h-1)u +\la x + \la xu = 0$. Therefore,  $y\in J_{\la,\mathrm{right}}$ and $g'\in G_\la$.

We have
\begin{equation}\label{xixi}
    \xitl(g') = \theta(h)\varepsilon^{\la((h-1)u+x+xu)}=\theta(h)\varepsilon^{\la(x)}=\xitl(g).
\end{equation}
Denote  $\La(g) = \{s\in G|~ s^{-1}gs\in G_\la\}$.

Let us prove that $\La(g') = \La(g)$. Really, $s^{-1}g's = 1+ (s^{-1}gs-1)s^{-1}as$.
As see above that if  $ s^{-1}gs\in G_\la$, then $s^{-1}g's\in G_\la$. This proves  $\La(g)\subset\La(g')$. The contrary can be proved similarly. The formula  (\ref{xixi}) implies \begin{equation}
\xitl(s^{-1}g's) = \xitl(s^{-1}g s)
\end{equation}
for any  $s\in \La(g)$.
Hence
\begin{equation}\label{chichi}
\chitl(g') = \frac{1}{|G_\la|}\sum_{s\in \La(g')} \xitl(s^{-1}g's)= \frac{1}{|G_\la|}\sum_{s\in \La(g)} \xitl(s^{-1}gs) = \chitl(g).
\end{equation}
{\bf Item 2}. The characters are constant on the classes of conjugate elements. Then  $\chitl(g_0gg_0^{-1})=\chitl(g)$ for any $g_0\in G$. Applying  (\ref{RRR}) and Item  1, we conclude
$\chitl(R_\tau(g))=\chitl(g)$ for any  $\tau\in\tG$ и $g\in G_\la$.

If a $\tG$-orbit has empty intersection with  $G_\la$, then $\chitl$ is zero on it.$\Box$\\
\Prop\Num\label{lalap}.
If $\la,~\la'\in J^*_e$ belong to a common  $\tGe$-orbit, then $\chitl=\chitlp$.\\
\Proof.  The statement follows from the following item 1 and 2.\\
{\bf Item 1}. Consider the case  $\la'=a\la$, where $a=1+u\in N_e$. Since  $ J_{a\la,\mathrm{right}}
=J_{\la,\mathrm{right}}$, we see $G_{a\la}=G_\la$. For any  $g=h+x\in G_\la$, we have
$$ \xitla(g) =   \theta(h)\varepsilon^{a\la(x)} = \theta(h)\varepsilon^{\la(x+xu)}= \theta(h)\varepsilon^{\la(x)}=\xitl(g),$$
where $g=h+x$, ~ $h\in H(e)$ and $x\in J_{\la,\mathrm{right}}$.
Hence, $\chitl=\chitla$ for any  $a\in N_e$.\\
{\bf Item 2}. Let $\la'= g_0\la g_0^{-1}$ where $g_0 = ta_0$, ~$t\in H$, $a_0\in N_e$. Then  $\Jlpr=g_0\Jlr g_0^{-1}$ and $g_0hg_0^{-1}=h$ for every $h\in H(e)$.
Therefore,  $G_{\la'}=g_0G_\la g_0^{-1}$ and
$$\xitlp(g) = \theta(h)\varepsilon^{g_0\la g_0^{-1}(x)} = \theta(h)\varepsilon^{\la(g_0^{-1}xg_0)} = \xitl(g_0^{-1}gg_0)$$
for any  $g\in G_{\la'}$. Then the following induced representations are equivalent and  $\chitl=\chitlp$. $\Box$\\
\Prop\Num\label{thethep}.
If $\theta\ne\theta'$, the  characters  $\chitl$ and $\chitpl$ are disjoint.\\
\Proof. Simplify notations:  $\xi=\xitl$ and $\xi'=\xitpl$. It follows from the Intertwining Number Theorem  \cite[теорема 44.5]{CR} that the characters  $\chitl$ and $\chitpl$ are disjoint if and only if  for any  $s\in G$ there exists  $g\in G_\la$ such that $sgs^{-1}\in G_\la$ and $\xi(sgs^{-1})\ne \xi'(g)$.

It is sufficient to prove that for any pair  $(a,t)$, where $a\in N$ and $t\in H$, there exist
$g\in G_\la$ such that  $aga^{-1}$ and $tgt^{-1}$ belong to  $G_\la$, ~  and
\begin{equation}\label{at}
    \xi(aga^{-1})\ne \xi'(tg t^{-1}).
\end{equation}\\
We shall call this property  \emph{the disjoint property} of the pair  $(a,t)$.\\
\emph{Remark}.
\begin{itemize}
 \item If  $a=a_1b$, where $a_1\in N$, ~ $b\in N_{\la,\mathrm{right}}$, and  $tbt^{-1}\in N_{\la,\mathrm{right}}$, then the disjoint property holds for  the pair  $(a,t)$ if and only if it holds for the pair  $(a_1,t)$.
\item  If  $a=ba_1$ where  $a_1\in N$, ~ $b\in N_{\la,\mathrm{right}}$, then   the disjoint property holds for  the pair $(a,t)$ if and only if it holds for the pair $(a_1,t)$.
\end{itemize}
{\bf Item  1}.  Let $a\in N$ and $t\in H$. We reduce the proof of  (\ref{at}) to the case $a=1+ev$, where  $v\in J$.
Let us show that there exists $w\in J$ such that  $a(1+e'w)=1+ev$, where  $v\in J$,~ $e'=1-e$. Indeed, let $a=1+u$,~ $u\in J$. For any  $v\in J$, we obtain $$ a(1+e'w) =  (1+eu+e'u)(1+e'w) = 1+ eu(1+e'w) + e'u + (1+e'u)e'w.$$
It suffices to take $e'w = -(1+e'u)^{-1}e'u$.

 As $\la e'=0$, we obtain  $e'w\in J_{\la,\mathrm{right}}$ and $1+e'w\in N_{\la,\mathrm{right}}$. Moreover,
$t(1+e'w)t^{-1}= 1+e' twt^{-1}\in N_{\la,\mathrm{right}}$. Applying the above Remark, we conclude that it is sufficient to prove the disjoint property in the case  $a=1+ev$, where  $v\in J$.\\
{\bf Item 2}. Let  $t\in H$ and  $a=1+ev=1+eve+eve'$, where $v\in J$.
We shall prove the disjoint property  for the case $eve'\in J_{\la,\mathrm{right}}$.

Then  $1-eve'\in N_{\la,\mathrm{right}}$ and   $(1-eve')a=(1-eve')(1+eve+eve') = 1+eve$.
Applying the above Remark, we conclude that it is sufficient to prove the disjoint property in the case  $a=1+eve$, where  $v\in J$.

Let  $g=h_0$ be en arbitrary element of the subgroup  $H(e)\subset G_\la$.
Then  $h_0e=eh_0=e$ and
$$ah_0a^{-1}= (1+eve)h_0(1+eve)^{-1} = h_0(1+eve)(1+eve)^{-1} =h_0\in G_\la.$$
Since $th_0=h_0t$, we have $th_0t^{-1} = h_0\in G_\la$.

  By assumption, $\theta\ne\theta'$; there  exists  $h_0\in H(e)$ such that $\theta(h_0)\ne\theta'(h_0)$. Substituting  $g=h_0$ in (\ref{at}), we obtain
 $$\xi(ah_0a^{-1}) = \theta(h_0) \ne \theta'(h_0) = \xi'(th_0t^{-1}).$$
 {\bf Item 3}. Let $a=1+x$, where $x=eve+eve'\in J$ and $eve'\notin J_{\la,\mathrm{right}}$.
Then $$\la(xe'J) =
\la((eve+eve')e'J)=\la((eve')e'J)=\la(eve'(eJ+e'J))=\la(eve'J)\ne 0.$$
Consider the chain of ideals  $e'J\supset e'J^2\supset\ldots\supset e'J^k=\{0\}.$
The exists a number  $i$ that  $\la(xe'J^i)\ne 0$ and $\la(xe'J^{i+1})= 0$.
 It follows, that there exists $y\in e'J^i$ such that  $\la(xy)\ne 0$ and $\la(xyJ) = 0$. The last equality means that  $xy\in J_{\la,\mathrm{right}}$. The element $y$ also belongs $ J_{\la,\mathrm{right}}$, since   $\la e'=0$.

Take $g=1+cy\in N_{\la,\mathrm{right}}$, where $c\in \Fq$. We have
$$ aga^{-1} =(1+x)(1+cy)(1+x)^{-1} =1+(1+x)cy(1+x)^{-1} \in 1+cy + cxy + yJ + xyJ$$
Therefore, $aga^{-1}\in N_{\la,\mathrm{right}}\subset G_\la.$
Notice that $\la(y)=0$, since  $\la(e'J)=0$. Also $\la(yJ) = \la(xyJ)=0$, since
 $y,~xy\in J_{\la,\mathrm{right}}$.
 We obtain
$\xi(aga^{-1}) = \varepsilon^{c\la(xy)}$.

As
$y\in e'J$, the element $tyt^{-1}$ also lies in $e'J$ and  $$tgt^{-1} = 1+tyt^{-1}\in N_{\la,\mathrm{right}}\subset G_\la.$$
We have   $\xi'(tgt^{-1}) =\varepsilon^{c\la(tyt^{-1})}=1$.

Since $\la(xy)\ne 0$, there exists  $c\in \Fq$ such that
$\varepsilon^{c\la(xy)}\ne 1$. This proves  (\ref{at}). $\Box$
\\
\Prop\Num\label{dva}. Let $\la,~\la'\in J_e^*$
and  $\theta$,~ $\theta'$ are linear characters of the group $H(e)$. Then\\
1) ~ $\chitl = \chitplp$ if and only if  $\la$ and $\la'$ belongs to a common  $\tGe$-orbit and  $\theta=\theta'$;\\
2)~ if  $\la$ and  $\la'$ belong to different  $\tGe$-orbits or  $\theta\ne\theta'$, then  the characters  $\chitl$ and  $\chitplp$ are disjoint.\\
\Proof. Recall the construction of supercharacters for the algebra group $N=1+J$ (see \cite{DI}). For any  $\mu\in J$, the supercharacter $\chi_\mu$ is an induced character from the character  $$\xi_\la(1+x)=\varepsilon^{\mu(x)}$$ of the subgroup $N_{\la,\mathrm{right}}$ to the group  $N$. Two supercharacters are equal if and only if $\mu$  and $\mu'$ are in the same  $N\times N$-orbit.  Otherwise the characters  $\mu$  and  $\mu'$ are disjoint  ~\cite[Theorem 5.5]{DI}.

Applying the formula for the  restriction of  induced representation  on the subgroup (see \cite[теорема 44.2]{CR}, \cite[proposition 22]{Serr}),  we obtain
\begin{equation}\label{res}
   \mathrm{Res}_N (\chitl) = \sum_{t\in H} \chi_{Ad_t\la}.
   \end{equation}
{\bf Case 1}. Let $\la$ and $\la'$ belongs to  a common  $\tGe$-orbit. Then  $\chitl=\chitlp$ (see proposition  \ref{lalap}).
 If $\theta=\theta'$, then $\chitplp =\chitlp=\chitl$. If  $\theta\ne\theta'$, then the  characters $\chitl=\chitlp$ and  $\chitplp$ are disjoint ( see proposition  \ref{thethep}). \\
{\bf Case 2}. Let $\la$ and $\la'$ belong to different  $\tGe$-orbits. Then  $\la$ and $\la'$ belong to different  $\tG$-orbits (see lemma \ref{Oc*}). By the formula (\ref{res}), the restrictions of  $\chitl$ and $\chitplp$ on the subgroup  $N$ are disjoint. Then the characters
$\chitl$ and $\chitplp$ are disjoint. $\Box$

Recall that   $\Ax$ stands for the set of triples  $\al = (e,\theta, \omega^*)$, where $e$ is an idempotent in  $kH$, ~ $\theta$ is a linear character of the group $H(e)$,  and
      $\omega^*$ is a  regular  $\tG_e$-orbit in  $J^*_e$.
For any  $\al\in\Ax$ denote by  $\chi_\al$ the supercharacter $\chitl$, where  $ \la\in\omega^*$.\\
 \Theorem\Num\label{main}. The systems of superchracters   $ \{\chi_\al|~ \al\in \Ax\}$ and superclasses  $\{\Kc_\beta|~ \beta\in \Bx\}$  form a supercharacter theory for the group    $G$.\\
 \Proof.  As we see above, the superchracters are pairwise disjoint. The superclasses form a partition of the group  $G$ (see Theorem  \ref{supcl}). Since the number of  regular orbits in  $J_e$ equals to the number of  regular orbits in  $J_e^*$ (see proposition \ref{ereg}), the number of supercharacters equals to the number of superclasses and equals to
\begin{equation}\label{AxAx}
|\Ax| = |\Bx| = \sum n_E(J_e)+|H(e)|,
\end{equation}
 where  $e$ runs through the set of all idempotents of $kH$  for whom   the set  of regular orbits in $J_e$ is nonempty.  This proves  S1). The supercharacters are constant on superclasses (see proposition  \ref{supersuper}); this proves  S2). Finally, $\{1\}$ is a superclass  $K(g)$ for $g=1$, this proves  S3). $\Box$

 \section{Examples}

{\bf Example 1}. The algebra group  $G=1+J$, where $J$ is a associative finite dimensional  algebra over the finite field  $k=\Fq$. Here $H=\{1\}$.  The constructed supercharacter theory coincides with the supercharacter theory of  P.Diaconis and  I.M.Isaacs  ~\cite{DI}. The superclass of an arbitrary element  $1+x$ is a set  $1+\omega$, where  $\omega $ is a left-right $G\times G$-orbit of  $x\in J$. The system of supercharacters consists of the characters $\chi_\la=\ind(\xi_\la, G_\la, G)$, where  $\la$ is a representative of left-right orbit in  $J^*$, ~ $G_\la$ is the stabiliser of of $\la$ with respect to the right  $G$-action,  and $\xi_\la(1+x)=\varepsilon^{\la(x)}$.\\
{\bf Example 2}. $G=\left\{ \left(\begin{array}{cc} a&b\\0&1\end{array}\right):~~ a,b\in k, ~a\ne 0\right\}$. The group  $G$ is a group of triangular type  $G=H+J$, where  $H=\left\{\left(\begin{array}{cc}a&0\\0&1\end{array}\right)\right\}$ and $J=\left\{\left(\begin{array}{cc}0&b\\0&0\end{array}\right)\right\}$.

 Choose a primitive root  $\zeta$ of degree  $q-1$ of  $1$ in the field  $k=\Fq$.
 The subgroup  $H$ is cyclic  $H=\{1,g,\ldots, g^{q-2}\}$, where $ g=\left(\begin{array}{cc}\zeta&0\\0&1\end{array}\right)$. Then the group algebra  $kH$ decomposes with respect to the system of primitive idempotents
$kH=ke_0\oplus ke_1\oplus\ldots ke_{q-2}$,  where
$$ e_i=\frac{1}{q-1}\sum_{j=0}^{q-2} \zeta^{-ij}g^j.$$

 Let $T: kH\to \mathrm{Mat}(2,k)$ be a representation of the algebra $kH$ that extends the  identical representation of the subgroup  $H$.
 Then $T(e_0)= \left(\begin{array}{cc}0&0\\0&1\end{array}\right)$, ~  $T(e_1)= \left(\begin{array}{cc}1&0\\0&0\end{array}\right)$, ~ $T(e_i)=0$ for all  $i>1$.

Let  $e$ be an idempotent in  $kH$, and $\mathrm{supp}(e)$ is a system of its primitive components.
If  $\mathrm{supp}(e)$ contains some primitive idempotent   $e_i$ with  $i>1$, then $J_e$ and $J_e^*$ has no  regular elements  (because  $J_e=J_{e-e_i}$).  Therefore,  these  idempotents are omitted  in the construction of  systems of parameters $\Ax$ and $\Bx$.

 If $e=e_0$ and  $e=e_1$, them easy to see that $J_e=0$, and again  $J_e$ has no  regular elements (because $J_e=J_f $ for $f=0$).
It remains to consider the cases  $e=0$ and  $e=e_0+e_1$.
If  $e=0$, then $J_e=0$, ~$H(e)=H$.  The supercharacters with  $e=0$  are linear characters of the group  $G$.

If $e=e_0+e_1$, then  $T(e) = 1$, ~ $J_e=J$, ~ $ H(e)=\{1\}$. A linear form  $\la\in J^*$ is regular if  it is nonzero. The set of nonzero elements of  $J^*$ form a  single $\rho(\tG)$-orbit. The supercharacter  $\chi_\la$ is induced by linear character of the group  $1+J$.

The system of supercharacters coincides with the system of irreducible characters, and the system of superclasses with the classes of conjugate elements.\\
{\bf Example 3}.
 Let $A$ be an associative finite dimensional algebra with unity over the finite field   $\Fq$ of $q$ elements \cite[\S 6.6]{Pi}. By definition, the  algebra $A$ is reduced if the its factor algebra with respect to the radical    $J=J(A)$ is a direct sum of division algebras. According to Wedderburn's theorem  \cite[\S 13.6]{Pi}, any division algebra over a finite field is commutative. Then the algebra  $A/J$ is commutative.
There exists a semisimple subalgebra $S$ such that
$A = S\oplus J$ (see \cite[\S 11.6]{Pi}). In our case  $S$ is commutative.
The group  $G=A^*$  on the invertible elements of  $A$ is a finite group of triangular type    $G=H+J$, where  $H=S^*$. If  $A$ is the algebra of triangular matrices, then $G=\Tn$ is the triangular group.

Consider the representation  $T: kH\to S$ that is identical on $H$. The image  $T(e_i)$ of each primitive idempotent  $e_i$ of  $kH$ is either zero, or is a primitive idempotent in $S$.
Similar to the previous example,  if   $e_i\in\mathrm{supp}(e)$ and  $T(e_i)=0$, then $J_e$ has no  regular elements. These idempotents are omitted  in construction of the systems of parameters $\Ax$ and $\Bx$.

The other idempotents  $e$ in $kH$ are in one-to-one correspondence with idempotents in  $S$.
Constructing the systems of parameters  $\Ax$ and $\Bx$,  one can consider that the idempotents are elements of  $S$.  Under this identification we see that the constructed in this paper supercharacter theory coincides in this example with the one of the paper  \cite{P1}.

\section{The analog of Kirillov's formula}

Recall the A.A.Kirillov formula for the character of irreducible representation of an unipotent group \begin{equation}\label{kir}
  \chi_\la (1+x) = \sum_{\mu \in\Omega^*} \varepsilon^{\mu(x)},
\end{equation}
where $\mu$ runs through the coadjoint orbit  $\Omega^*$ of $\la\in J^*$, and  $\chi_\la$
  is the irreducible representation  associated with $\Omega^*$. Remark that this formula is true if the characteristic of the field is sufficiently great \cite{Kir}, \cite{Ka}. The paper  $\cite{IsKar}$ implies that  the formula \ref{kir} is not correct for  the group  $\UTn$, ~ $n\ge 13$ if  $\mathrm{char}(\Fq)=2$.

In the paper  \cite{DI}, the analog of   Kirillov's formula for the supercharacters of the algebra groups was obtained
\begin{equation}\label{kirdi}
  \chi_\la (1+x) = \frac{1}{n(\la)}\sum_{\mu \in G\la G} \varepsilon^{\mu(x)},
\end{equation}
where $ \chi_\la$  is the supercharacter  associated with  $\la\in J^*$,
 and  $n(\la)$ is the number of the right  $G$-orbits in  $G\la G$.
In the same paper  \cite{DI}, there  proved the other version of  (\ref{kirdi}) in the form
\begin{equation}\label{kirdisec}
  \chi_\la (1+x) = \frac{|\la G|}{|GxG|}\sum_{y\in GxG} \varepsilon^{\la(y)}.
\end{equation}

Our goal is to prove the analog of there  A.A.Kirillov's formula for the supercharacters of the finite groups of  triangular type.
For this we need some lemmas.

 Let $h\in H$ and  $h-1\in kH$ associate with some idempotent  $f'\in kH$.  Then  $hf=fh=f$  for the idempotent $f=1-f'$. Denote $J_{11}=fJf=J_f$,~
$J_{12}=fJf'$,~ $J_{21}=f'Jf$,~ $J_{22}=f'Jf'$. We have the Pierce decomposition $$J=J_{11}\oplus J_{12}\oplus J_{21}\oplus J_{22}.$$
\Lemma\Num\label{sfirst}. Each element  $s=t+w$,  ~ $t\in H$,~ $w \in J$ uniquely  decomposes into a product  $s = (1+u)(t+v)$, where
$u=u_{12}\in J_{12}$  and  $v=v_{11}+v_{22}+v_{21}\in J_{11}\oplus J_{22}\oplus J_{21}$.\\
\Proof. From  $(1+u_{12})(t+v)=t+u_{12}t+v+u_{12}v$ we obtain
$w_{11}=u_{12}t+u_{12}v_{22}$,~ $w_{11}=v_{11}+u_{12}v_{21}$, ~ $w_{22}=v_{22}$, ~ $w_{21}=v_{21}$.
Knowing  $\{ w_{ij}\}$, we have  $v_{21}=w_{21}$,~ $v_{22}=w_{22}$,~ $u_{12}=w_{12}(t+w_{22})^{-1}$, и $v_{11}=w_{11}-u_{12}w_{21}$. $\Box$\\
\Lemma\Num\label{ssecond}. Let   $g=h+y$, where $h\in H$ and $y\in J$. The idempopent  $f$ is defined  via  $h$ as above.  Then\\
 1) there exist  $u=u_{12}\in J_{12}$  and  $z=z_{11}+z_{21}+z_{22}\in J_{11}\oplus J_{22}\oplus J_{21}$ such that  $h+y = (1+u)(h+z)(1+u)^{-1}$; the elements $u$ and  $z$ are uniquely defined by $g$.
\\
2) Given $\la\in J_f^*$,  then  $ y\in \Jlr$ if and only if  $u$ and $z$ belongs to $\Jlr$.\\
\Proof. \\
{ \bf Item 1}. Notice that  $(1+u_{12})^{-1} = 1-u_{12}$. We shall prove  that for any  $g=h+y$  there exists  a unique   $u=u_{12}\in J_{12}$ obeying  $$(1-u_{12})(h+y)(1+u_{12})= h+z,~~   z_{12}=0.$$

By direct calculations, we have
\begin{equation}\label{pmmm}\left\{\begin{array}{rll}
  z_{12}& =&y_{12} + u_{12}((1-h)-y_{22}-y_{21}u_{12}) +  y_{11}u_{12}, \\
  z_{11} &= &y_{11}-u_{12}y_{21}, \\
  z_{22}& = &y_{22}+y_{21}u_{12}, \\
  z_{21} &= &y_{21}.
\end{array}\right.
 \end{equation}
 Consider the map  $\phi: J_{12}\to J_{12}$ defined as  $$\phi(u_{12}) =
 u_{12}((1-h)-y_{22}-y_{21}u_{12}) +  y_{11}u_{12}.$$
The first equation of  (\ref{pmmm}) can be rewritten in the form  $ z_{12} =y_{12} + \phi(u_{12})$.

 Show that the map  $\phi$ is injective.
 Suppose that  $\phi(u_{12})=\phi(u_{12}')$. The equation  $\phi(u_{12})-\phi(u_{12}')=0$ implies
 \begin{equation}\label{llll}
   v(1-h)=vy_{22}-y_{11}v+vy_{21}u_{12}+u_{12}'y_{21}v,
 \end{equation}
 where $v= u_{12}-u_{12}'\in J_{12}$.

The nilpotent algebra $J$ has the natural filtration  $J\supset J^2\supset \ldots \supset J^M=\{0\}$.
 If $v\ne 0$, then  $v$ belongs to some  $J^m$ and does not belong to $J^{m+1}$. The equality  (\ref{llll}) implies   $v(1-h)\in J^{m+1}$.

By the other hand, the operator of right multiplication by  $1-h$ is invertible in $J_{12}$.
Therefore, the element т $v(1-h)$,  as well as  $v$, belongs to  $ J^m$ and does not belong to $J^{m+1}$. A contradiction.  Hence $u_{12}=u_{12}'$.

So, the map  $\phi$ is injective. Since the set  $J_{12}$ is finite, the map is $\phi$ a bijection.
There exists a unique  $u\in J_{12}$ such that  $z_{12}=y_{12} + \phi(u_{12})=0$.\\
{\bf Item 2}. By condition  $\la\in J_f$, then  $f\la=\la f=\la$. The element   $y$ lies in $\Jlr$ if and only if all its components  $y_{ij}$ are lying in $\Jlr$.

  Let $y\in \Jlr$. From  (\ref{llll}) we  obtain  $z_{21}, z_{22}\in \Jlr$. Take  $z_{12}=0$ in the first row of   (\ref{llll}). Then
 $$   u_{12}((1-h)-y_{22}-y_{21}u_{12}) + y_{12} + y_{11}u_{12}=0.$$
Therefore  $$   u_{12}((1-h)-y_{22}-y_{21}u_{12}) \in \Jlr.$$
The operator of right multiplication by  $1-h$ in invertible in $J_{12}$.
Hence, the operator of right multiplication by   $(1-h)-y_{22}-y_{21}u_{12}$  is also invertible in $J_{12}$.
This implies  $u_{12}\in \Jlr$. By the second row of   (\ref{pmmm}), we conclude $z_{11}\in \Jlr$.
So $u_{12}$ and  $z$ belongs to  $\Jlr$. The contrary assertion is proved similarly. $\Box$

Let  $\al = (e,\theta, \omega^*)\in \Ax$, where
 $e$ is an idempotent in   $kH$, ~ $\theta $  is a lonear character (one-dimensional representation) of the subgroup  $ H(e)$,  and  $\omega^*$ is a regular  $\tG_e$-orbit in $J^*_e$.
  Let $\chi_\al$ be the corresponding supercharacter.

  It follows from the theorem  \ref{superclass} each superclass has  the element $g=h+x$, where $hx=xh=x$. Since the supercharacters are constant on the superclasses, it is sufficient  to calculate values of supercharacters on the elements of these type.
 Notice that  the equality  $hx=xh=x$ is equivalent to  $f'x=xf'=0$,  that is   $x\in J_f$. Two elements  $g=h+x$ and  $g'=h+x'$, where $x,x'\in J_f$, belong to a common superclass if and only if $x$ and $x'$ belong to a common  $\rho(\tG_f)$-orbit (see theorem  \ref{superclass}). \\
\Lemma\Num\label{sthird}. If $h\notin H(e)$, then $\chi_\al(g)=0$ for any $g=h+x$,~ $x\in J_f$.\\
\Proof.  Choose  $\la\in\omega^*$. The supercharacter  $\chi_\al$ is induced from the linear character  $\xi_\la$ of the subgroup  $G_\la=H(e)+\Jlr$. Hence
\begin{equation}\label{chial}
\chi_\al(g) = \frac{1}{|G_\la|}\sum_{s\in G} \dot{\xi}_\la(sgs^{-1}).
\end{equation}
The element  $sgs^{-1}$ has the form  $h+y$, ~ $y\in J$. If $h\notin H(e)$, then $\chi_\al(g)=0$. $\Box$

Suppose that  $h\in H(e)$. Then  $e<f$ and, therefore,   $J_e^*$ is a subspace of  $J_f^*$.
In partiсular   $\omega^*\in J_f^*$. There exists a unique $\rho^*(\tG_f)$-orbit $\Omega^*$ in $J_f^*$ such that its intersection with  $J_e^*$ coincides with  $\omega^*$.
Denote by   $\dot{\theta}(h)$  the function on  $H\to\Cb$ that annihilates out off the subset  $H(e)$, and equal to   $\theta(h)$ on $H(e)$.\\
   \Theorem \Num\label{kirillov}.
The value of the supercharacter  $\chi_\al$ on the element  $g=h+x$, where  $hx=xh=x$, is calculated by as follows
  \begin{equation}\label{superkir}
\chi_\al(g) = \frac{|H_e|\cdot\dot{\theta}(h)}{n(\Omega^*)} \sum_{\mu\in \Omega^*} \varepsilon^{\mu(x)},
\end{equation}
where $n(\Omega^*)$ is the number of right  $N_f$-orbits in $\Omega^*$.\\
\Proof. By Lemma  \ref{sthird} it is sufficient to prove the  assertion for the case  $h\in H(e)$.
As above, we choose an arbitrary element  $\la\in\omega^*$.

Apply the formula  (\ref{chial}). The element  $sgs^{-1}$ is presented in the form  $h+y$, ~$y\in J$.
 As in Lemma  \ref{sfirst} we decompose   $s\in G$ into a product  $s=s_1s_2$, where $s_1=1+u_{12}$ and $s_2=t+v$, ~ $v=v_{11}+v_{22}+v_{21}$.
Hence  $s_2gs_2^{-1}=h+z$, where $z\in J_{11}\oplus J_{22}\oplus J_{21}$. The Lemma  \ref{ssecond} implies the representation  $h+y=s_1(h+z)s_1^{-1}$ is unique,  moreover  $ h+y\in G_\la$ (equivalent to  $y\in \Jlr$) if and only if
 $u_{12}\in\Jlr$ and   $ h+z\in G_\la$ (equivalent to  $z\in \Jlr$). We obtain
 $$ \chi_\al(g) = \frac{1}{|G_\la|}\cdot |(\Jlr)_{12}|\cdot\sum \dot{\xi}(s_2gs_2^{-1}),$$
where  $s_2$ runs through all elements of the form  $t+v$, ~ $t\in H$, ~$v\in J_{11}\oplus J_{22}\oplus J_{21}$.
Similar to Lemma  \ref{sfirst}, one can prove that  $s_2$  is uniquely decomposes into a  product $s_2 =\sigma_1\sigma_2\sigma$,  where $\sigma_1=1+u_{21}$,~ $\sigma_2=1+u_{22}$,~ $\sigma=t+u_{11}$. Notice that since $f'J\subset e'J\subset \Jlr$, the elements  $u_{21}$ and  $u_{22}$ belong to $\Jlr$.
Arguing as in Lemma  \ref{ssecond},  we verify that   $s_2 gs_2^{-1}$ belongs to  $G_\la$ if and only if $\sigma g\sigma^{-1}$ belongs to  $G_\la$.
Easy to see that  $\Jlr= (\Jlr)_{12} \oplus (\Jlr)_{11}\oplus f'J$; this implies  $$|\Jlr| = |(\Jlr)_{12}|\cdot |(\Jlr)_{11}|\cdot |f'J|.$$
We have
$$
 \chi_\al(g) = \frac{1}{|G_\la|}\cdot \frac{|\Jlr|}{|(\Jlr)_{11}|}
  \cdot \sum_{\sigma\in H+J_f}
  \dot{\xi}(\sigma g\sigma^{-1}).
$$
Taking into account  $|G_\la|=|H(e)|\cdot|\Jlr|$, we obtain
\begin{equation}\label{pppp}
 \chi_\al(g) =
\frac{\theta(h)}{|H(e)|\cdot|(\Jlr)_{11}|}
\sum_{t\in H}\sum
 \varepsilon^{t^{-1}\la t(a x a^{-1})},
\end{equation}
where in the second sum $a\in N_f$ runs through all elements obeying   $ axa^{-1}\in   J_{t^{-1}\la t,\mathrm{right}}$.

 Take $\nu=t^{-1}\la t$. The second sum in (\ref{pppp}) have the form
$$V(x)=\sum_{a\in N_f,~~ \nu axa^{-1}=0} \varepsilon^{\nu(a x a^{-1})}.$$
 Hence
$$ V(x)= \sum_{a\in N_f,~ \nu ax=0 } \varepsilon^{\nu a(x)} =
|(\Jlr)_{11}|\cdot\sum_{\mu\in \nu N_f,~ \mu x=0} \varepsilon^{\mu(x)}.$$

For any  $a_1\in N_f$ the right stabilizer of  $a_1\nu$ in $J$ coincides with the right stabilizer of  $\nu$.
That is  $J_{a_1\nu a,\mathrm{right}} = J_{\nu a,\mathrm{right}}$. Since  $x\in J_{\nu a,\mathrm{right}}$, we have   $\nu a(xJ)=0$,  and therefore
$a_1\nu a(x)= \nu a(xa_1)=\nu a(x)$.
Hence
$$\sum_{ a\in N_f,~ a_1\nu ax=0 } \varepsilon^{a_1\nu a(x)}=\sum_{ a\in N_f,~ \nu ax=0 } \varepsilon^{\nu a(x)} =V.$$
This implies that
\begin{equation}\label{pppm}
   V = \frac{|(\Jlr)_{11}|}{n(\la)}\sum_{\mu\in N_f\nu N_f,~ \mu x=0} \varepsilon^{\mu(x)},
 \end{equation}
where  $n(\la)$ is, as in  (\ref{kirdi}), the number if right  $N_f$-orbits in left-right orbit  $N_f\la N_f$.

Notice that for any nonzero linear form  $\la\in J^*$ the sum  $\sum_{u\in J}
\varepsilon^{\la(u)}$ equals to zero. If  $\mu x\ne 0$, then  $$\sum_{a_1=1+u_1,~ u_1 \in J_f}\varepsilon^{a_1\mu(x)} = \varepsilon^{\mu(x)} \sum_{u_1 \in J_f}\varepsilon^{\mu x(u_1)} = 0.$$
We can adopt the  requirement  $\mu x=0$ in  (\ref{pppm}).

 Substituting  (\ref{pppm})  into (\ref{pppp}), we obtain
\begin{multline}
\chi_\al(g) =\frac{\theta(h)}{|H(e)|\cdot n(\la)}\sum_{t\in H} \sum_{\mu\in N_f\la N_f} \varepsilon^{t^{-1}\mu t(x)}=\\
\frac{\theta(h)}{|H(e)|\cdot n(\la)}\cdot |\mathrm{Stab}(N_f\la N_f, H)|\sum_{\mu\in \Omega^*_f} \varepsilon^{\mu (x)}.
  \end{multline}
  The number of elements in the stabilizer  $\mathrm{Stab}(N_f\la N_f, H)$ equals to  $\frac{|H|}{n_1(\la)}$, where  $n_1(\la)$ is the number of left-right  $N_f\times N_f$-orbits in  $\tG_f$-orbit  $\Omega^*$.
  The equality $n_1(\la)n(\la)=n(\Omega^*)$ implies   $$
\chi_\al(g) =\frac{|H|\cdot \theta(h)}{|H(e)|\cdot n(\Omega^*)}\sum_{\mu\in \Omega^*}\varepsilon^{\mu(x)} $$
This proves the formula  (\ref{superkir}). $\Box$

In the next statement, we shall prove the version of the formula  (\ref{kirdisec}) for the finite groups of the triangular type.
\\
\Theorem\Num. The value of the supercharacter  $\chi_\al$ on the element $g=h+x$, where  $hx=xh=x$, is calculated by the formula:
\begin{equation}\chi_\al(g) = \frac{|H_e|\cdot |\la N_f|\cdot\dot{\theta}(h)}{|\rho(\tG_f)(x)|} \sum_{y\in \rho(\tG_f)(x) } \varepsilon^{\la(y)}.
\end{equation}
\Proof. Since $n(\Omega^*)=\frac{|\omega^*|}{|\la N_f|}$, the formula  (\ref{superkir}) have the form
 $$\chi_\al(g) = \frac{|H_e|\cdot |\la N_f|\cdot\dot{\theta}(h)}{|\Omega^*|} \sum_{\mu\in \Omega^*} \varepsilon^{\mu(x)}= \frac{|H_e|\cdot |\la N_f|\cdot\dot{\theta}(h)}{|\Omega^*|}\cdot\frac{|\Omega^*|}{|\tG_f|} \sum_{\tau\in \tG_f} \varepsilon^{\rho^*_\tau\la(x)} =$$
$$ \frac{|H_e|\cdot |\la N_f|\cdot\dot{\theta}(h)}{|\tG_f|} \sum_{\tau\in \tG_f} \varepsilon^{\la(\rho_\tau x)} =
 \frac{|H_e|\cdot |\la N_f|\cdot\dot{\theta}(h)}{|\tG_f|}\cdot\frac{|\tG_f|}{|\rho(\tG_f)(x)|}  \sum_{y\in \rho(\tG_f)(x) } \varepsilon^{\la(y)} = $$ $$ \frac{|H_e|\cdot |\la N_f|\cdot\dot{\theta}(h)}{|\rho(\tG_f)(x)|} \sum_{y\in \rho(\tG_f)(x) } \varepsilon^{\la(y)}.~~\Box
 $$

\end{document}